\documentclass[a4paper,12pt]{article}
\usepackage{latexsym,amsmath}
\usepackage{indentfirst}
\usepackage{graphicx}
\begin{document}

\begin{center}
{\large \bf{ FRACTIONAL DYNAMICAL SYSTEMS ON FRACTIONAL LEIBNIZ ALGEBROIDS }}\\
\end{center}
\begin{center}
{\bf  Gheorghe IVAN, Mihai IVAN and Dumitru OPRI\c S }\\[0.2cm]
\end{center}
\def \b{\Box}

\textit{Dedicated to Acad. Prof. Dr. Radu Miron at his 80 th
anniversary }\\[0.2cm]

\begin{center}
Gheorghe Ivan, West University of Timi\c soara\\
Departments of Mathematics\\
4, B-dul V. P{\^a}rvan, 300223, Timi\c{s}oara, Romania\\
E-mail: ivan @ math.uvt.ro\\
\end{center}

\begin{center}
Mihai Ivan, West University of Timi\c soara\\
Seminarul de Geometrie - Topologie\\
4, B-dul V. P{\^a}rvan, 300223, Timi\c{s}oara, Romania\\
E-mail: mihai31ro @ yahoo.com\\
\end{center}

\begin{center}
Dumitru  Opri\c s, West University of Timi\c soara\\
Departments of Mathematics\\
4, B-dul V. P{\^a}rvan, 300223, Timi\c{s}oara, Romania\\
E-mail: miticaopris @ yahoo.com\\
\end{center}

 {\bf Abstract.} In this paper we consider the fractional tangent bundle on a differentiable manifold. A fractional Leibniz structure on an
  algebroid is defined.The fractional dynamical system on a fractional Leibniz algebroid is defined and it is discussed. Some illustrative examples are presented.{\footnote{ \hspace*{0.7cm}\textit{Key words and phrases}. Fractional
derivatives, fractional tangent bundle, fractional Leibniz
algebroid, fractional differential equations}

\section { Introduction}

The theory of derivative of noninteger order goes back to Leibniz,
Liouville, Riemann, Grunwald and Letnikov. Derivatives of
fractional order have found many applications in recent studies in
mechanics, physics, economics, medicine.Classes of fractional
differentiable systems have studied in [10], [4].

In the first section the fractional tangent bundle to a
differentiable manifold is defined, using the method of Radu
Miron's from [8].
 In this paper the fractional dynamical
systems on fractional Leibniz algebroids are presented. The
associated geometrical objects have an geometric character. Also,
some examples for fractional dynamical systems of this type are
given.

\section{Fractional tangent bundle on a manifold }

Let $ f: [a,b] \rightarrow \textbf{R} $ and $ \alpha \in
\textbf{R}, \alpha
> 0 $. The \textit{ Riemann - Liouville fractional derivative at
to left of $ a $ }, respectively \textit{ at to right of $ b $ }
is the function $ f \rightarrow _{a}D_{t}^{\alpha}f $ resp. $ f
\rightarrow _{t}D_{b}^{\alpha}f $, where:

\begin{equation}
\left\{\begin{array}{lcl}
_{a}D_{t}^{\alpha}f(t)&=&\frac{1}{\Gamma(m - \alpha )}(\frac{d}{dt})^{m}\int_{a}^{t}{(t-s)^{m-\alpha -1}}(f(s) - f(a))ds\\[0.2cm]
_{t}D_{b}^{\alpha}f(t)&=&\frac{1}{\Gamma (m -
\alpha)}(-\frac{d}{dt})^{m}\int_{t}^{b} {(t-s)^{m-\alpha -1
}}(f(s) - f(b))ds,
\end{array}\right.
\label {1}
\end{equation}
where $ m \in \textbf {N}^{\ast} $ such that $ m-1 \leq \alpha
\leq m , \Gamma $ is the Euler gamma function and $(
\frac{d}{dt})^{m}=\frac{d}{dt}\circ \frac{d}{dt} \circ ...\circ
\frac{d}{dt}. $

 We will denote sometimes $ D_{t}^{\alpha} =
_{a}D_{t}^{\alpha}.$

The following proposition holds.

\textbf{Proposition 2.1 }([3]) \textit{ (i) If}
$\lim_{n\rightarrow \infty} \alpha_{n} =p \in \textbf {N}^{\ast}
$,\textit{ then:}
\begin{equation}
 \lim_{n\rightarrow
\infty}(_{a}D_{t}^{\alpha_{n}}f(t))= D_{t}^{p}f(t), ~~~~~
\lim_{n\rightarrow
\infty}(_{t}D_{b}^{\alpha_{n}}f(t))= D_{t}^{p}f(t).\\[0.2cm]
. \label {2}
\end{equation}

\textit{ (ii) If $ f(t) = c, (\forall)  t \in [a,b], ~
D_{t}^{\alpha}f(t) = 0 . $}

\textit{ (iii) If $ f_{1}(t) = t^{\gamma}, (\forall)  t \in
[a,b],$ then $ ~ D_{t}^{\alpha}f_{1}(t) =
\frac{\Gamma(1+\gamma)}{\Gamma (1+\gamma - \alpha)} t^{\gamma -
\alpha }. $ }

\textit{ (iv) If $ f_{1}, f_{2} $ are analytical functions on $
(a,b ), $  then:}

\begin{equation}
 D_{t}^{\alpha}(f_{1} f_{2})(t)  =
\sum_{k=0}^{\infty}\left (\begin{array}{c} \alpha\\
k
\end{array}\right )D_{t}^{\alpha - k}f_{1}(t)
(\frac{d}{dt})^{k}f_{2}(t).\\[0.2cm]
\label {3}
\end{equation}

\textit{ (v) If $ f: [a,b] \rightarrow\textbf {R} $ is analytical
function on $ (a,b), $ and $ 0\in (a,b) $ then:}

\begin{equation}
 f(t)  =  \sum_{h=0}^{\infty}
 {E_{\alpha,h}(t) D_{t}^{\alpha  h} f(t)}|_{t=0},\\[0.2cm]
\label {4}
\end{equation}
\textit{where $ E_{\alpha,h} $ is the Mittag - Leffler 's
function:}
\begin{equation}
 E_{\alpha,h}(t)  =  \sum_{h=0}^{\infty}\frac{1}{\Gamma(1+\alpha h)}t^{\alpha  h
} .\\[0.2cm]
\label {5}
\end{equation}
$\hfill \b$

Let $ \alpha \in \textbf {R}, \alpha
> 0 $ and $ M $ be a manifold of dimension $ n $ and $
U $ a local chart on $ M $. We say that the curves $ c_{1}, c_{2}
: I \rightarrow M, 0\in I, c_{1}(0)=c_{2}(0)\in M $ have \textit{
fractional contact $\alpha $  in $ x_{0} $,} if for all $ f\in
C^{\infty}(U), x_{0}\in U, $ the following relation holds:
\begin{equation}
 D_{t}^{\alpha} (f\circ c_{1})|_{t=0} = D_{t}^{\alpha} (f\circ c_{2})|_{t=0}
 .\label {6}
\end{equation}

 The set of equivalences classes $ ([c]_{x_{0}}^{\alpha}) $ is
 called the \textit{ fractional tangent space in $ x_{0} $ } and
 it is denoted by $ T_{x_{0}}^{\alpha}(U) $.

 Let  $ T^{\alpha}(M)= \bigcup_{x_{0}\in M} T_{x_{0}}^{\alpha}(U) $
 and the projection $ \pi^{\alpha} : T^{\alpha}(M) \rightarrow M $ given by
$ \pi^{\alpha}([c]_{x_{0}}^{\alpha})=x_{0} .$

 On $ T^{\alpha}(M)$ there exists a differentiable structure and
 we can prove that $ ( T^{\alpha}(M), \pi^{\alpha}, M ) $ is a
 differentiable bundle.

In a system of local coordinates on $ M ,$ if $ x_{0}\in U $ and $
c:I\rightarrow M $ is a curve given by $ x^{i}=x^{i}(t), (\forall
) t\in I ,$   the class $ ( [c]_{x_{0}}^{\alpha}) $ is given by:
\begin{equation}
x^{i}(t)= x^{i}(0) + \frac{1}{\Gamma(1+\alpha)}t^{\alpha}
D_{t}^{\alpha}x^{i}(t)|_{t=0},~~~ t\in (\varepsilon, \varepsilon
). \label {7}
\end{equation}

On the open set $ (\pi^{\alpha})^{-1} (U)\in T^{\alpha}(M), $ the
local coordinates of the element $ ( [c]_{x_{0}}^{\alpha}) $ are $
(x^{i}, y^{i(\alpha)} ), $ where:
\begin{equation}
x^{i}= x^{i}(0), ~~~~~ y^{i(\alpha )} = \frac{1}{\Gamma(1+\alpha)}
D_{t}^{\alpha}x^{i}(t),  i=\overline{1,n}. \label {8}
\end{equation}

\textbf{ Proposition 2.2} ([1], [2] ) \textit{ Let  $ U,
\overline{U} $ be two local charts on $ M $ such that $ U \cap
\overline{U}\neq \emptyset $  and
\begin{equation}
\overline{x}^{i}= \overline{x}^{i}( x^{1}, x^{2}, ..., x^{n} ),
~~~ det(\frac{\partial \overline{x}^{i}}{\partial x^{j}})\neq 0,
~~~ i=\overline{1,n} \label {9}
\end{equation}
the coordinate transformations. The coordinate transformations on
$ (\pi^{\alpha})^{-1}(U\cap \overline{U}) $ are given by:
\begin{equation}
\overline{x}^{i}= \overline{x}^{i}(x^{1}, x^{2},..., x^{n}), ~~~
\overline{y}^{i(\alpha )} = {\overset{\alpha}{J}}_{j}^{i}(x,
\overline{x}) y^{j(\alpha)}, \label {10}
\end{equation}
 where:
 \begin{equation}
\overset {\alpha}{J}_{j}^{i}(x, \overline{x})=
\frac{1}{\Gamma(1+\alpha )} D_{\overline{x}^{j}}^{\alpha}
(x^{i})^{\alpha}\label {11}
\end{equation}
 and  $ D_{x^{i}}^{\alpha} $ is defined by:
\begin{equation}
\begin{array}{l}
D_{x^{i}}^{\alpha}f(x)=\frac{1}{\Gamma(1-\alpha
)}\frac{\partial}{\partial x^{i}}\cdot\\[0.4cm]
~~~~~~~~~~~\cdot \int_{a^{i}}^{x^{i}}
\frac{f(x^{1},...,x^{i-1},s,x^{i+1},...,x^{n})-
f(x^{1},...,x^{i-1},a^{i},x^{i+1},...,x^{n}) }{(x^{i}-
s)^{\alpha}} ds ,
\end{array}
\label {12}
\end{equation}
with $ U_{ab} = \{ x\in U, a^{i} \leq x^{i} \leq b^{i},
i=\overline{1,n} \}\subseteq U. $ }$\hfill\b$

Let $ \mathcal{D}^{1}(U) $  the module of $ 1 $ - forms defined on
$ U \subseteq M. $ Using the  fractional exterior derivative $
d^{\alpha} : C^{\infty}(U) \rightarrow \mathcal{D}^{1}(U), ~ f
\rightarrow d^{\alpha}(f) $ ( see [2] ), where $ d^{\alpha}(f) $
is given by:
 \begin{equation}
d^ {\alpha}(f)= d(x^{i})^{\alpha} D_{x^{i}}^{\alpha}(f)\label {13}
\end{equation}
 follows:
 \begin{equation}
 \left\{\begin{array}{lcl}
d(x^{i})^{\alpha}& = &
{\overset{\alpha}{J}}_{j}^{i}(x,\overline{x})
d(\overline{x}^{j})^{\alpha}\\[0.2cm]
D_{x^{i}}^{\alpha} & = &
{\overset{\alpha}{J}}_{i}^{j}(\overline{x},x)D_{\overline{x}^{j}}^{\alpha}\\[0.2cm]
{\overset{\alpha}{J}}_{k}^{i}(x,\overline{x}){\overset{\alpha}{J}}_{j}^{k}(\overline{x},x)&=&
\delta_{j}^{i}.
\end{array}\right.
\label{14}
\end{equation}

We denote by $ \mathcal{X}^{\alpha}(U) $ the module of fractional
vector fields generated by the operators $ \{ D_{x^{i}}^{\alpha},
~i=\overline{1,n} \} .$ A fractional vector field $
{\overset{\alpha}{X}}\in \mathcal{X}^{\alpha}(U) $  has the
following form:
\begin{equation}
{\overset{\alpha}{X}} = {\overset{\alpha}{X}}^{i}
D_{x^{i}}^{\alpha}, ~~~~~ {\overset{\alpha}{X}}^{i}\in
C^{\infty}(U), i=\overline{1,n}, \label{15}
\end{equation}
which for a change of local charts, the correspondent components
satisfies the relations:
\begin{equation}
  {\overset{\alpha}{\overline {X}}}^{i} =
  {\overset{\alpha}{J}}_{j}^{i}(x,\overline{x}){\overset{\alpha}{X}}^{j},~~~ i,j=\overline{1,n}.\label{16}
\end{equation}

The fractional differentiable equations associated to fractional
vector field $ {\overset{\alpha}{X}} $ is:
\begin{equation}
D_{t}^{\alpha} x^{i}(t) =
  {\overset{\alpha}{X}}^{i}(x(t)),~~~ i=\overline{1,n} \label{17}
\end{equation}
or equivalently ( using the notation (8)):
\begin{equation}
\Gamma ( 1+\alpha) y^{i(\alpha)} (t)=
  {\overset{\alpha}{X}}^{i}(x(t)),~~~ i=\overline{1,n}. \label{18}
\end{equation}

The fractional differential equations $ (17) $ with initial
conditions have solutions, see [3]. Examples of fractional
differentiable equations on $ \textbf {R} $ can be find in [4].

\section{ Fractional Leibniz dynamical systems }

Let the module $ \mathcal{X}^{\alpha}(U) $ of fractional vector
fields generated by the operators $ \{ D_{x^{i}}^{\alpha},
~i=\overline{1,n} \} $ and the module $ \mathcal{D}^{\alpha}(U) $
 generated by the $ 1 - $ forms $ \{ d (x^{i})^{\alpha},
~i=\overline{1,n} \}. $
 Applying the Proposition 2.2 it follows:
 \begin{equation}
 (d(x^{i})^{\alpha})(
 D_{x^{j}}^{\alpha})=D_{x^{j}}^{\alpha}(x^{i})^{\alpha} =
 \Gamma(1+\alpha) \delta_{j}^{i}. \label{19}
 \end{equation}

 If $ {\overset{\alpha}{X}}\in \mathcal{X}^{\alpha}(U) $ and $ {\overset{\alpha}{\omega}}\in \mathcal{D}^{\alpha}(U) $
 such that $
 {\overset{\alpha}{\omega}}={\overset{\alpha}{\omega}}_{i}d(x^{i})^{\alpha} ,$ then $
{\overset{\alpha}{\omega}}({\overset{\alpha}{X}})=\Gamma(1+\alpha){\overset{\alpha}{X}}^{i}{\overset{\alpha}{\omega}}_{i}.$

Let be a  fractional $ 2- $ contravariant tensor field $
{\overset{\alpha}{B}}\in \mathcal{X}^{\alpha}(U)\times
\mathcal{X}^{\alpha}(U) $ and $ d^{\alpha} f, d^{\alpha}g \in
\mathcal{D}^{\alpha}(U) $ defined by $(13) .$

The bilinear map $ [\cdot,\cdot]^{\alpha} : C^{\infty}(M) \times
C^{\infty}(M) \rightarrow C^{\infty}(M)$  defined by:
\begin{equation}
[f,g]^{\alpha} = B( d^{\alpha}f,d^{\alpha}g ), ~~~ (\forall )
f,g\in C^{\infty}(M), \label{20}
\end{equation}
is called the \textit{ fractional Leibniz bracket}.

If $ {\overset{\alpha}{B}}= {\overset{\alpha}{B}}^{ij}
D_{x^{i}}^{\alpha}\otimes D_{x^{j}}^{\alpha} ,$ from $(20) $
follows:
\begin{equation}
[f,g]^{\alpha} = {\overset{\alpha}{B}}^{ij}\cdot
D_{x^{i}}^{\alpha}f \cdot D_{x^{j}}^{\alpha} g .\label{21}
\end{equation}

Since
\begin{equation}
 D_{x^{i}}^{\alpha}(f h)(x)  =
\sum_{k=0}^{\infty}\left (\begin{array}{c} \alpha\\
k
\end{array}\right )(D_{x^{i}}^{\alpha - k}f(x))
(\frac{\partial}{\partial x^{i}})^{k}h(x),\\[0.2cm]
\label {22}
\end{equation}
it follows
\begin{equation}
 [f h, g ]^{\alpha}  =
\sum_{k=0}^{\infty}\left (\begin{array}{c} \alpha\\
k
\end{array}\right )\cdot {\overset{\alpha}{B}}^{ij}(D_{x^{i}}^{\alpha -
k}f)\cdot
(\frac{\partial}{\partial x^{i}})^{k}h \cdot D_{x^{j}}^{\alpha} g,\\[0.2cm]
\label {23}
\end{equation}

Similarly, one obtain
\begin{equation}
 [f , g h ]^{\alpha}  =
\sum_{k=0}^{\infty}\left (\begin{array}{c} \alpha\\
k
\end{array}\right )\cdot
{\overset{\alpha}{B}}^{ij}(D_{x^{i}}^{\alpha}f)\cdot
D_{x^{j}}^{\alpha - k}(g)\cdot
(\frac{\partial}{\partial x^{i}})^{k}h .\\[0.2cm]
\label {24}
\end{equation}

The pair $ ( M, [\cdot, \cdot]^{\alpha} ) $ is called \textit{
fractional Leibniz manifold}. If the bracket $ [\cdot,
\cdot]^{\alpha} $ is skew-symmetric, that is  $ [f, g]^{\alpha} =
- [g , f]^{\alpha} $ for all $ f,g\in C^{\infty}(M) $ we say that
$ ( M, [\cdot, \cdot]^{\alpha} ) $ is a \textit{ fractional almost
Poisson manifold}. If $ \alpha \rightarrow 1 ,$ then one obtain
the concept from  [6].

For $ h\in C^{\infty}(M) ,$ the fractional vector field $
{\overset{\alpha}{X}}_{h} $ defined by
\begin{equation}
{\overset{\alpha}{X}}_{h}(f)= [f,h]^{\alpha}, ~ (\forall) f\in
C^{\infty}(M), \label{25}
\end{equation}
is called the \textit{ fractional Leibniz vector field associated
to} $ h. $ The fractional dynamical system associated to $
{\overset{\alpha}{X}}_{h} $ is called the \textit{ fractional
Leibniz dynamical system}.

 If $ (x^{i}), i=\overline{1,n} $ is a system of local coordinates
 on $ M $, then the fractional
Leibniz dynamical system is given by
\begin{equation}
D_{t}^{\alpha} x^{i}(t)= [x^{i}(t),h(t)]^{\alpha}, ~~~\textrm{
where}~~~ [x^{i},h]^{\alpha} = {\overset{\alpha}{B}}^{ij}\cdot
D_{x^{j}}^{\alpha}h  . \label{26}
\end{equation}

\textbf{Example 3.1.} Let the constant fractional $ 2-$
contravariant tensor $ {\overset{\alpha}{g}} = (
{\overset{\alpha}{g}}^{ij}) $ defined on $ \textbf{R}^{3} $ by
\begin{equation}
{\overset{\alpha}{g}} = \left ( \begin{array}{ccc}
s_{1}\gamma_{1}& 0 & 0\\
0 & s_{2} \gamma_{2} & 0\\
0 & 0 & s_{3} \gamma_{3}
\end{array}\right ), \label{27}
\end{equation}
where $ s_{1}, s_{2}, s_{3} \in \{-1,1\} $ and $ \gamma_{1},
\gamma_{2}, \gamma_{3} $ satisfies the relation $ \gamma_{1}+
\gamma_{2}+ \gamma_{3}= 0. $

For  $ h = x^{1}x^{2}x^{3}, $ the associated fractional Leibniz
dynamical system is
\begin{equation}
\left\{\begin{array}{lclcl} D_{t}^{\alpha} x^{1} & = &
s_{1}\gamma_{1} D_{x^{1}}^{\alpha}(h) &=&
\frac{\Gamma(2)}{\Gamma(2-\alpha)}s_{1}\gamma_{1}x^{2} x^{3}
(x^{1})^{1-\alpha}\\[0.2cm]

D_{t}^{\alpha} x^{2} & = & s_{2}\gamma_{2} D_{x^{2}}^{\alpha}(h)
&=& \frac{\Gamma(2)}{\Gamma(2-\alpha)}s_{2}\gamma_{2}
x^{1} x^{3}(x^{2})^{1-\alpha}\\[0.2cm]

D_{t}^{\alpha} x^{3} & = & s_{3}\gamma_{3} D_{x^{3}}^{\alpha}(h)
&=& \frac{\Gamma(2)}{\Gamma(2-\alpha)}s_{3}\gamma_{3}x^{1} x^{2}
(x^{3})^{1-\alpha}. \label{28}
\end{array}\right.
\end{equation}

If $ \alpha \rightarrow 1 ,$ it follows the system $(7) $  in
[10].

For  ${\overset{\alpha}{h}} = (x^{1})^{\alpha} (x^{2})^{\alpha }
(x^{3})^{\alpha}, $ the associated fractional Leibniz dynamical
system is
\begin{equation}
\left\{\begin{array}{lcl}
 D_{t}^{\alpha} x^{1} & = & \Gamma(1+\alpha)
s_{1}\gamma_{1} x^{2} x^{3}\\[0.2cm]
D_{t}^{\alpha} x^{2} & = & \Gamma(1+\alpha)
s_{2}\gamma_{2} x^{1} x^{3}\\[0.2cm]
D_{t}^{\alpha} x^{3} & = & \Gamma(1+\alpha) s_{3}\gamma_{3} x^{1}
x^{2} . \label{29}
\end{array}\right.
\end{equation}
$\hfill\b$

Let $ {\overset{\alpha}{P}} $ be a  skew-symmetric fractional $ 2-
$ contravariant tensor field and a non-degenerate symmetric
fractional $ 2-$ contravariant tensor field $
{\overset{\alpha}{g}} $ on  the manifold $ M$. We define the
bracket $ [\cdot, \cdot]^{\alpha} : C^{\infty}(M) \times
C^{\infty}(M) \rightarrow C^{\infty}(M)$ by
\begin{equation}
[f,h]^{\alpha} = {\overset{\alpha}{P}}( d^{\alpha}f,d^{\alpha}h
)+{\overset{\alpha}{g}}( d^{\alpha}f,d^{\alpha}h ) , ~~~ (\forall
) f,h\in C^{\infty}(M). \label{30}
\end{equation}

The $ 4-$ tuple $ ( M,{\overset{\alpha}{P}},
{\overset{\alpha}{g}}, [\cdot, \cdot]^{\alpha} ) $ is called
\textit{ fractional almost metric manifold}.

The fractional dynamical system associated to $ h\in C^{\infty}(M)
$ is
\begin{equation}
D_{t}^{\alpha} x^{i}(t)= [x^{i}(t),h(t)]^{\alpha}, ~~~\textrm{
where}~~~ [x^{i},h]^{\alpha} = {\overset{\alpha}{P}}^{ij}
D_{x^{j}}^{\alpha}h + {\overset{\alpha}{g}}^{ij}
D_{x^{j}}^{\alpha}h . \label{31}
\end{equation}

\textbf{Example 3.2.} Let be the fractional $ 2 $ - contravariant
tensors fields $ {\overset{\alpha}{P}} =
({\overset{\alpha}{P}}^{ij} )  , {\overset{\alpha}{g}} = (
{\overset{\alpha}{g}}^{ij} ) $ on $ {\bf R}^{3} $ and the function
$ h\in C^{\infty}({\bf R}^{3} ) $ given by :

$${\overset{\alpha}{P}} = \left ( \begin{array}{ccc}
0 & x^{3} & - x^{2}\\
-x^{3} & 0 & x^{1}\\
x^{2} & - x^{1}& 0\\
\end{array}\right ),$$\\[-0.9cm]

$${\overset{\alpha}{g}} =\left ( \begin{array}{ccc}
-a_{2}(x^{2})^{2}- a_{3}(x^{3})^{2} & a_{1} a_{2} x^{1} x^{2} & a_{1} a_{3} x^{1} x^{3}\\
a_{1} a_{2} x^{1} x^{2} & -a_{1}(x^{1})^{2}- a_{3}(x^{3})^{2} & a_{2} a_{3} x^{2} x^{3} \\
a_{1} a_{3} x^{1} x^{3} & a_{2} a_{3} x^{2} x^{3} & -a_{1}(x^{1})^{2}- a_{2}(x^{2})^{2} \\
\end{array}\right ),$$\\[-0.9cm]

$$ h = ( a_{1} + 1 )( x^{1})^{\alpha} + ( a_{2} +
1 )( x^{2})^{\alpha} + ( a_{3} + 1 )( x^{3})^{\alpha}.$$\\[-0.5cm]

Since $ D_{x^{1}}^{\alpha}h = (a_{1}+1) \Gamma(1+\alpha),
D_{x^{2}}^{\alpha}h = (a_{2}+1) \Gamma(1+\alpha),\\
D_{x^{3}}^{\alpha}h = (a_{3}+1) \Gamma(1+\alpha) ,$ the fractional
Leibniz dynamical system (31) associated to $ h $ is
$$ \left ( \begin{array} {c}
D_{t}^{\alpha} x^{1}\\
D_{t}^{\alpha} x^{2}\\
D_{t}^{\alpha} x^{3}\\
\end{array} \right ) = \Gamma(1+\alpha){\overset{\alpha}{P}}\left ( \begin{array} {c}
a_{1}+1\\
a_{2}+1\\
a_{3}+1\\
\end{array} \right )+ \Gamma(1+\alpha){\overset{\alpha}{g}}\left ( \begin{array} {c}
a_{1}+1\\
a_{2}+1\\
a_{3}+1\\
\end{array} \right
)=$$\\[-0.5cm]
$$=\Gamma(1+\alpha)({\overset{\alpha}{P}}+{\overset{\alpha}{g}})\left
( \begin{array} {c}
a_{1}+1\\
a_{2}+1\\
a_{3}+1\\
\end{array} \right
). \hfill\b $$

Let be two fractional $ 2 $ - contravariant tensors fields $
{\overset{\alpha}{P}}$ and $ {\overset{\alpha}{g}}$ on $ M.$
Define the bracket $ [\cdot, (\cdot, \cdot)] : C^{\infty}(M)\times
C^{\infty}(M)\times C^{\infty}(M)\rightarrow C^{\infty}(M)$ by:
\begin{equation}
[f,h]^{\alpha} = {\overset{\alpha}{P}}(
d^{\alpha}f,d^{\alpha}h_{1} )+{\overset{\alpha}{g}}(
d^{\alpha}f,d^{\alpha}h_{2} ) , ~~~ (\forall ) f,h_{1}, h_{2}\in
C^{\infty}(M). \label{32}
\end{equation}

The fractional vector field  ${\overset{\alpha}{X}}_{h_{1}h_{2}}$
defined by
\begin{equation}
 {\overset{\alpha}{X}}_{h_{1}h_{2}} =
[f, (h_{1},h_{2})],~~~ (\forall ) f\in C^{\infty}(M). \label{33}
\end{equation}
 is called the \textit{ fractional almost Leibniz vector field associate to the functions
 } $ h_{1}, h_{2} \in C^{\infty}(M).$ The dynamical system
 associated to ${\overset{\alpha}{X}}_{h_{1}h_{2}}$ is called the \textit{ fractional almost Leibniz dynamical
 system}.

 Locally, the fractional almost Leibniz dynamical
 system is given by:
\begin{equation}
D_{t}^{\alpha} x^{i}(t)= {\overset{\alpha}{P}}^{ij}
D_{x^{j}}^{\alpha}h_{1} + {\overset{\alpha}{g}}^{ij}
D_{x^{j}}^{\alpha}h_{2} . \label{34}
\end{equation}

\textbf{Example 2.3.} Let be the fractional $ 2 $ - contravariant
tensors fields $ {\overset{\alpha}{P}} =
({\overset{\alpha}{P}}^{ij} )  , {\overset{\alpha}{g}} = (
{\overset{\alpha}{g}}^{ij} ) $ on $ {\bf R}^{3} $ and the
functions $ h_{1}, h_{2}\in C^{\infty}({\bf R}^{3} ) $ given by :

$${\overset{\alpha}{P}} = \left ( \begin{array}{ccc}
0 & 1 & 0\\
- 1 & 0 & x^{1}\\
0 & - x^{1}& 0\\
\end{array}\right ),$$\\[-0.9cm]

$${\overset{\alpha}{g}} =\left ( \begin{array}{ccc}
0 & 0 & 0\\
0 & - (x^{3})^{2} & 0 \\
0 & 0 & -(x^{2})^{2} \\
\end{array}\right ),$$\\[-0.9cm]

$$ {\overset{\alpha}{h}}_{1} = ( x^{2})^{1+\alpha} + (
x^{3})^{1+\alpha},~~ {\overset{\alpha}{h}}_{2} = (
x^{1})^{1+\alpha} + ( x^{3})^{\alpha}.$$\\[-0.5cm]

Since\\[0.2cm]
$ D_{x^{1}}^{\alpha}{\overset{\alpha}{h}}_{1}=0, ~~~
D_{x^{2}}^{\alpha}{\overset{\alpha}{h}}_{1}=\Gamma(1+\alpha)
x^{2},~~~
D_{x^{3}}^{\alpha}{\overset{\alpha}{h}}_{1}=\Gamma(1+\alpha)
x^{3};$\\[0.2cm]
$ D_{x^{1}}^{\alpha}{\overset{\alpha}{h}}_{2}=\Gamma(1+\alpha)
x^{1}, ~~~ D_{x^{2}}^{\alpha}{\overset{\alpha}{h}}_{2}=0,~~~
D_{x^{3}}^{\alpha}{\overset{\alpha}{h}}_{2}=\Gamma(1+\alpha)
,$\\[0.2cm]
 the system (34) becomes:\\[0.2cm]
 $\left ( \begin{array}{c}
 D_{t}^{\alpha}x^{1}\\
D_{t}^{\alpha}x^{2}\\
D_{t}^{\alpha}x^{3}\\
\end{array}\right )= \left ( \begin{array}{ccc}
0 & 1 & 0\\
- 1 & 0 & x^{1}\\
0 & - x^{1}& 0\\
\end{array}\right )\left ( \begin{array}{c}
 0\\
\Gamma(1+\alpha)x^{2}\\
\Gamma(1+\alpha)x^{3}\\
\end{array}\right )+$\\[0.4cm]
$~~~~~~~~~~~~~~~~~~~~+\left ( \begin{array}{ccc}
0 & 0 & 0\\
0 & - (x^{3})^{2} & 0 \\
0 & 0 & -(x^{2})^{2} \\
\end{array}\right )\left ( \begin{array}{c}
 \Gamma(1+\alpha)x^{1}\\
0\\
\Gamma(1+\alpha)\\
\end{array}\right ) $ \\[0.2cm]

or equivalently

\begin{equation}
\left \{ \begin{array}{lcl}
D_{t}^{\alpha}x^{1} & = & \Gamma(1+\alpha) x^{2}\\
D_{t}^{\alpha}x^{2} & = & \Gamma(1+\alpha) x^{1}x^{3}\\
D_{t}^{\alpha}x^{3} & = & - \Gamma(1+\alpha) x^{1}x^{2}-
\Gamma(1+\alpha) (x^{2})^{2}. \label{35}
\end{array}\right.
\end{equation}

 The system (35) is called the \textit{ fractional
Maxwell- Bloch equations}.

If in (35) we take $ \alpha \rightarrow 1 , $ then one obtain the
Maxwell-Bloch equations.

\section{ Fractional Leibniz algebroids }

Let $ M $ be a smooth manifold of dimension $ n $ , let $ \pi : E
\rightarrow M $ be a vector bundle and $  \pi^{\ast} : E^{\ast}
\rightarrow M $ the dual vector bundle. By $ Sec(M,E) $ or $
Sec(\pi) $  we denote the sections of $ \pi $.

A \textit{ fractional Leibniz algebroid structure } on a vector
bundle $ \pi : E \rightarrow M $ is
 given by a bracket ( bilinear operation ) $ [ \cdot, \cdot ]^{\alpha} $ on the space of sections $ Sec(\pi) $ and two vector bundle morphisms
$ {\overset{\alpha}{\rho}}_{1}, {\overset{\alpha}{\rho}}_{2} : E
\rightarrow T^{\alpha}M $ ( called the \textit {left} and the
\textit { right fractional anchor} , respectively ) such that
\begin{equation}
\left \{ \begin{array}{l} [e_{a}, e_{b}]^{\alpha}
=C_{ab}^{c}e_{c}\\[0.2cm]
 [ f \sigma_{1}, g \sigma_{2} ]^{\alpha}
=f {\overset{\alpha}{\rho}}_{1}(\sigma_{1})(g)\sigma_{2} - g
{\overset{\alpha}{\rho}}_{2}(\sigma_{2})(f) \sigma_{1} + f g [
\sigma_{1}, \sigma_{2}]^{\alpha}\\
\end{array}\right.
 \label { 36}
\end{equation}
for all $ \sigma_{1}, \sigma_{2}\in Sec(\pi)$ and $ f, g \in
C^{\infty}(M).$

 A vector bundle  $ \pi : E \rightarrow M $ endowed with a fractional Leibniz algebroid structure $ ( [\cdot, \cdot ]^{\alpha}, {\overset{\alpha}{\rho}}_{1}, {\overset{\alpha}{\rho}}_{2} ) $ on $ E $ , is called {\it
  fractional Leibniz algebroid} over $ M $ and denoted by $(E, [\cdot, \cdot ]^{\alpha}, {\overset{\alpha}{\rho}}_{1}, {\overset{\alpha}{\rho}}_{2} ).$

 A fractional Leibniz algebroid with an antisymmetric bracket $ [\cdot, \cdot ]^{\alpha} $ ( in this case we have $ {\overset{\alpha}{\rho}}_{1}=- {\overset{\alpha}{\rho}}_{2} $ ) is called \textit {fractional pre - Lie algebroid}.

In a system of local coordinates the relation (36) reads:
\begin{equation}
[ \sigma_{1}^{a}e_{a},\sigma_{2}^{b}e_{b}]^{\alpha}
=\sigma_{1}^{a}{\overset{\alpha}{\rho}}_{1}(e_{a})(\sigma_{2}^{b})e_{b}-
\sigma_{2}^{a}{\overset{\alpha}{\rho}}_{2}(e_{a})(\sigma_{1}^{b})e_{b}+
\sigma_{1}^{a}\sigma_{2}^{b}C_{ab}^{c}e_{c}.\label{37}
\end{equation}

If $
{\overset{\alpha}{\rho}}_{1}(e_{a})={\overset{\alpha}{\rho}}_{1a}^{i}D_{x^{i}}^{\alpha},
~
{\overset{\alpha}{\rho}}_{2}(e_{b})={\overset{\alpha}{\rho}}_{2b}^{i}D_{x^{i}}^{\alpha}
$, from $(37)$ follows:
\begin{equation}
[ \sigma_{1}^{a}e_{a},\sigma_{2}^{b}e_{b}]^{\alpha}
=\sigma_{1}^{a}{\overset{\alpha}{\rho}}_{1a}^{i}(D_{x^{i}}^{\alpha}\sigma_{2}^{b})e_{b}-
\sigma_{2}^{a}{\overset{\alpha}{\rho}}_{2a}^{i}(D_{x^{i}}^{\alpha}\sigma_{1}^{b})e_{b}+
\sigma_{1}^{a}\sigma_{2}^{b}C_{ab}^{c}e_{c}.\label{38}
\end{equation}

In the following, we establish a correspondence between the
fractional Leibniz algebroid structures on the vector bundle $ \pi
: E \to M $ and the fractional $ 2 $- contravariant tensor fields
on bundle manifold $ E^{*} $ of the dual vector bundle $ \pi^{*} :
E^{*} \to M $.

For a given section $ \sigma \in Sec(\pi) ,$ we define the
function $ i_{E^{*}}\sigma $ on $ E^{*} $ by the relation :
\begin{equation}
i_{E^{\ast}}\sigma(a) = < \sigma(\pi^{\ast}(a)),a > ,~~\hbox { for
} ~~ a\in E^{\ast},  \label {39}
\end{equation}
where $ <\cdot, \cdot > $ is the canonical pairing between $ E $
and $ E^{\ast}. $
 If $ \sigma = \sigma^{a}e_{a} $ and $ a\in E^{\ast}$ has the
 coordinates $ ( x^{i}, \xi_{a} ), $ then:
 \begin{equation}
i_{E^{\ast}}\sigma(a) = \sigma^{a} \xi_{a}.  \label {40}
\end{equation}

Let $ {\overset{\alpha}{\Lambda}} $ be a fractional $ 2 $ -
contravariant tensor field on $ E^{\ast} $ and the bracket $
[\cdot, \cdot ]_{{\overset{\alpha}{\Lambda}}} $ of functions
defined by:
 \begin{equation}
 [ f, g ]_{{\overset{\alpha}{\Lambda}}^{\beta}} = {\overset{\alpha}{\Lambda}}^{\beta}( d^{\alpha \beta}f, d^{\alpha \beta}g ), ~(\forall) ~ f,g\in C^{\infty}(E^{\ast}),\label {41}
 \end{equation}
where
 \begin{equation}
 d^{\alpha \beta}f = d(x^{i})^{\alpha} D_{x^{i}}^{\alpha}f + d(\xi_{a})^{\beta}D_{\xi_{a}}^{\beta}f= d^{\alpha}(f)+d^{\beta}(f).\label {42}
 \end{equation}

 In the basis  $\{ D_{x^{i}}^{\alpha}, D_{\xi_{a}}^{\beta}\}, i=\overline{1,n},
 a=\overline{1,m} $ of the module $
 \mathcal{X}^{\alpha\beta}({\pi^{\ast}}^{-1}(U)),$ the components  $ \Lambda^{\alpha\beta} $ are given by:
\begin{equation}
\Lambda^{\alpha\beta}= A_{ab}D_{\xi_{a}}^{\beta}\otimes
D_{\xi_{b}}^{\beta} + A_{1a}^{i}D_{\xi_{a}}^{\beta}\otimes
D_{x^{i}}^{\alpha} + A_{2a}^{i}D_{x^{i}}^{\alpha}\otimes
D_{\xi_{a}}^{\beta}. \label{43}
\end{equation}

 For  a given fractional $ 2 $ - contravariant tensor field   $ {\overset{\alpha}{\Lambda}}^{\beta} $ on $ E^{\ast}, $ we say that $ {\overset{\alpha}{\Lambda}}^{\beta} $
 is \textit{ linear}, if for each pair $ ( \mu_{1}, \mu_{2} ) $ of sections of $ \pi^{\ast} $ , the function
 $ {\overset{\alpha}{\Lambda}}^{\beta} ( d (i_{E^{\ast}}\mu_{1})^{\beta}, d (i_{E^{\ast}}\mu_{2})^{\beta} ) $ defined on $ E^{\ast} $ is linear with respect the coordinates $ \xi_{a}.$

If $ \mu_{1}=\mu_{1}^{a}(x)e_{a}, \mu_{2}=\mu_{2}^{a}(x)e_{a},$
then $ d_{E^{\ast}}\mu_{1}=\mu_{1}^{a}(x)\xi_{a},
d_{E^{\ast}}\mu_{2}=\mu_{2}^{a}(x)\xi_{a}$ and
 $ {\overset{\alpha}{\Lambda}}^{\beta} ( d (i_{E^{\ast}}\mu_{1})^{\beta}, d (i_{E^{\ast}}\mu_{2})^{\beta} )=A_{ab}(x,\xi)
 (\mu_{1}^{c}(x))^{\alpha}(\mu_{2}^{a}(x))^{\alpha})D_{\xi_{a}}^{\beta}(\xi_{c})D_{\xi_{b}}^{\beta}(\xi_{a})^{\beta}=
\frac{1}{\Gamma(1+\alpha)^{2}}(\mu_{1}^{a}(x))^{\beta}(\mu_{2}^{b}(x))^{\beta}A_{ab}(x,\xi).
  $

 It follows that $ {\overset{\alpha}{\Lambda}}^{\beta}$ is linear
 if and only if $ A_{ab}(x,\xi) = C_{ab}^{c}(x)\xi_{c}.$

The fractional formulation of the Grabowski and Urbanski's Theorem
from [6], is the following.

 \textbf {Theorem 4.1.}\textit{ For every fractional Leibniz algebroid structure on $
\pi : E \to M $
 with the bracket $ [ \cdot, \cdot ]^{\alpha} $ and the fractional anchors $ {\overset{\alpha}{\rho}}_{1}, {\overset{\alpha}{\rho}}_{2} $ , there exists an unique
 fractional $ 2 $ - contravariant tensor field
$ {\overset{\alpha}{\Lambda}} $ on $ E^{\ast} $ such that the
following relations hold:}
\begin{equation}
\left \{ \begin{array}{lcl}
 i_{E^{\ast}}[ \sigma_{1}, \sigma_{2} ] &=&
[ (i_{E^{\ast}}\sigma_{1})^{\beta},
(i_{E^{\ast}}\sigma_{2})^{\beta}
]_{{\overset{\alpha}{\Lambda}}^{\beta}}\\[0.2cm]
\pi^{\ast}({\overset{\alpha}{\rho}}_{1}(\sigma)(f)) &=& [
(i_{E^{\ast}}\sigma)^{\beta} ,\pi^{\ast}
f]_{{\overset{\alpha}{\Lambda}}^{\beta}}\\[0.2cm]
\pi^{\ast}({\overset{\alpha}{\rho}}_{2}(\sigma)(f)) &=&
  [ \pi^{\ast} f, (i_{E^{\ast}}\sigma)^{\beta}
]_{{\overset{\alpha}{\Lambda}}^{\beta}},  \label {44}
\end{array}\right.
\end{equation}

\textit{ for all $ \sigma, \sigma_{1}, \sigma_{2} \in Sec(\pi)$
and $ f\in C^{\infty}(M). $}

\textit{ Conversely, every  fractional $ 2 $ - contravariant
linear tensor field $ {\overset{\alpha}{\Lambda}}^{\beta} $ on $
E^{\ast} $ defines a fractional Leibniz algebroid on $ E $ if the
relations ( 44 )  hold.} $\hfill\b$

Let $ ( x^{i} ), i=\overline{1,n} $ be a local coordinate system
on $U\subseteq M $ and let $ \{ e_{1}, \ldots, e_{m} \} $ be a
basis of local sections of $ E|_{U} $ ( $ dim ~M = n, dim~ E = n
+m $ ). We denote by $ \{ e^{1}, \ldots, e^{m}\} $ the dual basis
of local sections of $ E^{\ast}|_{U} $ and $ ( x^{i}, y^{a} ) $ (
resp., $ ( x^{i}, \xi_{a} ) $ ) the corresponding coordinate
system on $ E $ ( resp., $ E^{\ast} $ ).

Let $ {\overset{\alpha}{\Lambda}}^{\beta} $ given by (43). Using
(44),  it is easy to see that every linear fractional $ 2 $ -
contravariant tensor field $ {\overset{\alpha}{\Lambda}}^{\beta} $
on $ E^{\ast} $ has the form:
\begin{equation}
{\overset{\alpha}{\Lambda}}^{\beta} = C_{ab}^{d}\xi_{d}
D_{\xi_{a}}^{\beta}\otimes D_{\xi_{b}}^{\beta} +
{\overset{\alpha}{\rho}}_{1 a}^{i} D_{\xi_{a}}^{\beta}\otimes
D_{x^{i}}^{\alpha}  -
{\overset{\alpha}{\rho}}_{2a}^{i}D_{x^{i}}^{\alpha}\otimes
D_{\xi_{a}}^{\beta} , \label {45}
\end{equation}

where $ C_{ab}^{d}, {\overset{\alpha}{\rho}}_{1 a}^{i},
{\overset{\alpha}{\rho}}_{2 a}^{i}\in C^{\infty}(M) $  are
functions of $ x^{i} $.

The correspondence  between $ {\overset{\alpha}{\Lambda}}^{\beta}
$ and a fractional Leibniz algebroid structure is given by the
following relations :
\begin{equation}
[ e_{a}, e_{b}]^{\alpha} = C_{ab}^{d}e_{d}~,~
{\overset{\alpha}{\rho}}_{1}(e_{a}) = {\overset{\alpha}{\rho}}_{1
a}^{i} D_{x^{i}}^{\alpha}~,~ {\overset{\alpha}{\rho}}_{2}(e_{a})
={\overset{\alpha}{\rho}}_{2 a}^{i} D_{x^{i}}^{\alpha}.\label {46}
\end{equation}

We call a \textit{ fractional dynamical system on the fractional
Leibniz algebroid $ \pi : E\to M $ }, the fractional dynamical
system associated to vector field $
{\overset{\alpha}{X}}_{h}^{\beta} $ with $ h\in
C^{\infty}(E^{\ast}) $ given by:
\begin{equation}
{\overset{\alpha}{X}}_{h}^{\beta}(f) =
{\overset{\alpha}{\Lambda}}^{\beta}(d^{\alpha\beta}f,
d^{\alpha\beta}h),~\hbox { for all }~ f\in
C^{\infty}(E^{\ast}).\label {47}
\end{equation}

In a system of local coordinates $ ( x^{i}, \xi_{a} ) $ on $
E^{\ast} $, the dynamical system $ (47) $ is given by :

\begin{equation}
\left\{ \begin{array}{c}
 D_{t}^{\alpha} \xi_{a} = [ \xi_{a},
h]_{{\overset{\alpha}{\Lambda}}^{\beta}}\\[0.3cm]
D_{t}^{\alpha}x^{i} = [ x^{i},
h]_{{\overset{\alpha}{\Lambda}}^{\beta}}
\end{array}\right.
.\label {48}
\end{equation}
where
\begin{equation}
\left\{ \begin{array}{l}
[\xi_{a},h]_{{\overset{\alpha}{\Lambda}}^{\beta}} =
C_{ab}^{d}\xi_{d}D_{\xi_{b}}^{\beta} h +
{\overset{\alpha}{\rho}}_{1a}^{i}D_{x^{i}}^{\alpha} h\\[0.3cm]
 [ x^{i},h]_{{\overset{\alpha}{\Lambda}}^{\beta}} = -{\overset{\alpha}{\rho}}_{2a}^{i}D_{\xi_{a}}^{\beta} h
\end{array}\right.
.\label {49}
\end{equation}

If $\alpha \rightarrow 1, \beta \rightarrow 1 ,$ dynamical system
(48) was studied in [6].

If $\alpha \rightarrow 1$ dynamical system (48) has the form:
\begin{equation}
\left\{\begin{array}{l} \dot{x}^{i} = -
{\overset{\alpha}{\rho}}_{2a}^{i} D_ {\xi_{a}}^{\beta} h\\[0.3cm]
D_{t}^{\beta} \xi_{a}= C_{ab}^{d}\xi_{d}D_{\xi_{b}}^{\beta} h +
{\overset{\alpha}{\rho}}_{1a}^{i}\frac{\partial h}{\partial x^{i}}
\end{array}\right.
.\label {50}
\end{equation}

If $\beta \rightarrow 1$ dynamical system (48) has the form:
\begin{equation}
\left\{\begin{array}{lcl}
 D_{t}^{\alpha}x^{i}& =& -
{\overset{\alpha}{\rho}}_{2 a}^{i}\frac{\partial h}{\partial
\xi_{a}}\\[0.3cm]
 \dot{\xi}_{a}&=& C_{ab}^{d}\xi_{d}\frac{\partial
h}{\partial \xi_{b}} + {\overset{\alpha}{\rho}}_{1 a}^{i}
D_{x^{i}}^{\alpha}h
\end{array}\right.
.\label {51}
\end{equation}

If the fractional Leibniz algebroid is a fractional pre - Lie
algebroid ( that is, $ C_{ab}^{d} = - C_{ba}^{d} $ ), then the
 fractional dynamical system $ (48 )$ is given by  :
\begin{equation}
\left \{\begin{array}{lcl}
  D_{t}^{\beta} \xi_{a}& =&
C_{ab}^{d}\xi_{d} D_{\xi_{a}}^{\beta}h +
{\overset{\alpha}{\rho}}_{1 a}^{i}D_{ x^{i}}^{\alpha}h\\[0.3cm]
D_{t}^{\alpha} x^{i}& =& -{\overset{\alpha}{\rho}}_{1 a}^{i}D_{
\xi_{a}}^{\beta}h
\end{array}\right.
.\label {52}
\end{equation}

If the fractional Leibniz algebroid is a fractional symmetric
algebroid ( that is,$ C_{ab}^{d} =  C_{ba}^{d} $) , then the
fractional dynamical system $ (48)$ is given by:
\begin{equation}
\left\{\begin{array}{lcl}
  D_{t}^{\beta} \xi_{a}& =&
C_{ab}^{d}\xi_{d} D_{\xi_{a}}^{\beta}h +
{\overset{\alpha}{\rho}}_{1a}^{i}D_{x^{i}}^{\alpha}h\\[0.3cm]
D_{t}^{\alpha} x^{i}& =& {\overset{\alpha}{\rho}}_{1a}^{i}D_{
\xi_{a}}^{\beta}h
\end{array}\right.
.\label {53}
\end{equation}

\textbf {Example 4.1.} Let the vector bundle $ \pi : E = {\bf
R}^{3}\times {\bf R}^{3}\to {\bf R}^{3} $ and
 $ \pi^{*} : E^{*}={\bf R}^{3}\times ({\bf R}^{3})^{*}\to {\bf R}^{3} $ the dual vector bundle. We consider on $ E^{*} $ the fractional $ 2 $ - contravariant linear
 tensor field $ {\overset{\alpha}{\Lambda}} $  defined by the matrix $ P^{\beta} $, the fractional anchors $ {\overset{\alpha}{\rho}}_{1},{\overset{\alpha}{\rho}}_{2} $ and the function $ h $ given
 by:\\[0.2cm]
$P^{\beta} = \left ( \begin{array}{ccc}
0 & -\xi_{3}x^{3} & \xi_{2}x^{2}\\
\xi_{3}x^{3} & 0 & -\xi_{1}x^{1}\\
- \xi_{2}x^{2} & \xi_{1}x^{1} & 0\\
\end{array}\right ) ~~,~~{\overset{\alpha}{\rho}}_{1} =\left ( \begin{array}{ccc}
0 &- x^{3} & x^{2}\\
x^{3} & 0 & 0 \\
- x^{2} & 0 & 0 \\
\end{array}\right ),$\\[0.4cm]

${\overset{\alpha}{\rho}}_{2} =\left ( \begin{array}{ccc}
0 &-1 & 0\\
1 & 0 & -x^{1} \\
0 & x^{1} & 0  \\
\end{array}\right )~$ and  $~ h(x,\xi) =  (x^{2})^{\alpha}(\xi_{2})^{\beta} +
(x^{3})^{\alpha}(\xi_{3})^{\beta},  \alpha >0, \beta>0.$

Using the calculus formulas:\\[0.2cm]
$ D_{\xi_{a}}^{\beta}(\xi_{b})^{\gamma} =\delta_{b}^{a}
\xi_{a}^{\gamma-\beta}\frac{\Gamma(1+\gamma)}{\Gamma(1+\gamma-\beta)}
~~ ,~~D_{x^{i}}^{\alpha}(x^{j})^{\gamma}=\delta_{i}^{j}
(x^{i})^{\gamma-\alpha}\frac{\Gamma(1+\gamma)}{\Gamma(1+\gamma-\alpha)}$\\[0.2cm]
follows:
$$\left\{\begin{array}{ccc}
D_{\xi_{1}}^{\beta} h=0, & D_{\xi_{2}}^{\beta}
h=\Gamma(1+\beta)(x^{2})^{\alpha},& D_{\xi_{3}}^{\beta}
h=\Gamma(1+\beta)(x^{3})^{\alpha}\\
D_{x^{1}}^{\alpha} h=0, & D_{x^{2}}^{\alpha}
h=\Gamma(1+\alpha)(\xi_{2})^{\beta},& D_{x^{3}}^{\alpha}
h=\Gamma(1+\alpha)(x^{3})^{\beta}.\\
\end{array}\right.$$

The fractional dynamical system  ( 48 ) for the given elements,
has the following matrix form:\\[0.2cm]

$\left (\begin{array}{l}
D_{t}^{\beta}\xi_{1}\\
D_{t}^{\beta}\xi_{2}\\
D_{t}^{\beta}\xi_{3}\\
\end{array}\right) =  \Gamma(1+\beta)\left ( \begin{array}{ccc}
0 & -\xi_{3}x^{3} & \xi_{2}x^{2}\\
\xi_{3}x^{3} & 0 & -\xi_{1}x^{1}\\
- \xi_{2}x^{2} & \xi_{1}x^{1} & 0\\
\end{array}\right)\left (\begin{array}{c}
0\\
(x^{2})^{\alpha}\\
(x^{3})^{\alpha}\\
\end{array}\right)+$\\[0.4cm]
 $~~~~~~~~~~~~~~+\Gamma(1+\alpha)\left ( \begin{array}{ccc}
0 &- x^{3} & x^{2}\\
x^{3} & 0 & 0 \\
- x^{2} & 0 & 0 \\

\end{array}\right)\left (\begin{array}{c}
0\\
(\xi_{2})^{\beta}\\
(\xi_{3})^{\beta}
\end{array}\right)
,$\\[0.2cm]

$\left (\begin{array}{c}
D_{t}^{\alpha}x^{1}\\
D_{t}^{\alpha}x^{2}\\
D_{t}^{\alpha}x^{3}\\
\end{array}\right)=- \Gamma(1+\beta)\left ( \begin{array}{ccc}
0 & - 1 & 0\\
1 & 0 & -x^{1}\\
0 & x^{1} & 0\\
\end{array}\right)\left (\begin{array}{c}
0\\
(x^{2})^{\alpha}\\
(x^{3})^{\alpha}\\
\end{array}\right).$

From the above matrix equations  follows the fractional dynamical
system:
\begin{equation}
\left\{\begin{array}{lcl}
D_{t}^{\beta} \xi_{1} & = &
\Gamma(1+\beta)( -\xi_{3}(x^{2})^{\alpha}x^{3} +
\xi_{2}x^{2}(x^{3})^{\alpha})+\\
 &&+\Gamma(1+\alpha)(
-x^{3}(\xi_{2})^{\beta}+ x^{2}(\xi_{3})^{\beta})\\
D_{t}^{\beta} \xi_{2} & = & -\Gamma(1+\beta)
\xi_{1}(x^{3})^{\alpha}\\
D_{t}^{\beta} \xi_{3} & = & -\Gamma(1+\beta)
\xi_{1}(x^{2})^{\alpha}\\
D_{t}^{\beta} \xi_{3} & = & \Gamma(1+\beta)
\xi_{1}(x^{2})^{\alpha}\\
D_{t}^{\alpha}x^{1}& = &- \Gamma(1+\beta)(x^{2})^{\alpha}\\
D_{t}^{\alpha}x^{2}& = &- \Gamma(1+\beta)x^{1}(x^{3})^{\alpha}\\
D_{t}^{\alpha}x^{3}& = & \Gamma(1+\beta)x^{1}(x^{3})^{\alpha}\\
\end{array}\right.
.\label{54}
\end{equation}

The fractional dynamical system (54) is the
\textit{$(\alpha,\beta)- $ fractional dynamical system} associated
to fractional Maxwell-Bloch equations.$\hfill\b$

\textbf{Conclusion.} The numerical integration of the fractional
systems presented in this paper will be discussed in  future
papers.

\begin{center}
{\bf References}\\
\end{center}
\hspace*{0.7cm} [1]. I.D. Albu, M. Neam\c tu, D. Opri\c s,\textit{
The geometry of fractional osculator bundle of higher order and
applications.} Conference of Differential Geometry: Lagrange and
Hamilton spaces. Dedicated to Acad. Prof. dr.
Radu Miron at his 80-th anniversary. September 3 - 8, 2007, Ia\c si, Romania.\\
\hspace*{0.7cm} [2]. K. Cotrill- Shepherd, M. Naber, \textit{ Fractional differential forms }. J. Math. Phys., {\bf 42} (2001), 2203-2212.\\
\hspace*{0.7cm} [3]. J. Cresson, \textit{ Fractional embedding of
differential operators and Lagrangian systems .} Journal of Mathematical Physics, {\bf 38},Issue 3, 2007.\\
\hspace*{0.7cm} [4]. D. Deac, D. Opri\c s , \textit{ The geometry
of fractional osculator bundle of higher order on \textbf{R}.
Applications to mechanics and economics.}Conference of
Differential Geometry: Lagrange and Hamilton spaces. Dedicated to
Acad. Prof. dr.
Radu Miron at his 80-th anniversary. September 3 - 8, 2007, Ia\c si, Romania.\\
\hspace*{0.7cm} [5]. J. Grabowski, P. Urbanski, \textit{ Lie
algebroid and Poisson - Nijenhuis structures.} Rep. Math. Phys., {\bf 40}, 1997, 195 - 208.\\
\hspace*{0.7cm} [6].  Gh. Ivan, D. Opri\c s,\textit{ Dynamical
systems on Leibniz algebroids}.Differential Geometry - Dynamical systems (DGDS),\textbf{8} (2006). Geometry Balkan Press,127-137.\\
\hspace*{0.7cm} [7]. Gh. Jum\u arie, \textit{ Lagrangian mechanics
of fractional order, Hamilton- Jacobi fractional PDE and Taylor's
Series of nondifferentiable functions.} Chaos, Solitons and Fractals, 2007.\\
\hspace*{0.7cm} [8]. R. Miron , \textit{ The geometry of higher
order Lagrange spaces. Applications to Mechanics and Physics.} Kluwer Academic Publisher FTPH, {\bf 82},1997.\\
\hspace*{0.7cm} [9]. J.P. Ortega, V. Planas - Bielsa , \textit{ Dynamics on Leibniz manifolds .} Preprint, arXiv: math. DS/ 0309263, 2003.\\
\hspace*{0.7cm} [10]. V.E. Tarasov, \textit{ Fractional generalization of gradient and Hamiltonian systems.} J. Phys. A. Math. Gen., {\bf 38} ( 2005), 5929 - 5943.\\

 \end{document}